\documentclass[12pt,leqno]{article}

\def\int{\mathbb{Z}}

\def\vep{\varepsilon}
\def\Ueg{{\cal U}_{\varepsilon}(G)}
\def\Uegchi{{\cal U}_\chi(G)}
\def\Uel{{\cal U}_{\varepsilon}(L)}
\def\Uep{{\cal U}_{\varepsilon}(P)}

\def\Ueh{{\cal U}_{\varepsilon}(T)}
\def\Uehchi{{\cal U}_{\chi_T}(T)}
\def\Uelchi{{\cal U}_\chi(L)}
\def\Uelchil{{\cal U}_{\chi_L}(L)}
\def\Uegxi{{\cal U}_\xi(G)}
\def\Uepchi{{\cal U}_\chi(P)}
\def\Uepchil{{\cal U}_{\chi_L}(P)}
\def\mO{{\mathcal O}}
\def\Chi{\Xi}
\def\zog{Z_0(G)}
\def\zol{Z_0(L)}
\def\g{\mathfrak g}

\def\l{{\mathfrak l}}
\def\h{{\mathfrak h}}

\def\proof{{\bf Proof. }}

\def\pf{\proof}

\def\spec{{\rm Spec}}
\usepackage{times}
\usepackage{graphics}
\usepackage{amssymb}
\usepackage{amscd}
\usepackage{amsmath}
\input xy 
\xyoption{all}
\title{Induced conjugacy classes and induced $\Ueg$-modules}

\newtheorem{theorem}{Theorem}[section]
\newtheorem{conjecture}{Conjecture}
\newtheorem{lemma}[theorem]{Lemma}
\newtheorem{corollary}[theorem]{Corollary}
\newtheorem{proposition}[theorem]{Proposition}
\newtheorem{definition}[theorem]{Definition}
\newtheorem{remark}[theorem]{Remark}

\author{Giovanna Carnovale\\
Dipartimento di Matematica\\
via Trieste 63 - 35121 Padova - Italy\\
email: carnoval@math.unipd.it}

\begin{document}

\maketitle

\begin{abstract}
By work of De Concini, Kac and Procesi the irreducible representations of the non-restricted specialization of the quantized enveloping algebra of the Lie algebra $\g$ at the roots of unity are parametrized by the conjugacy classes of a group $G$ with ${\rm Lie}(G)=\g$. 
We show that there is a natural dimension preserving bijection between the sets of irreducible representations associated with conjugacy classes lying in the same Jordan class (decomposition class). We conjecture a relation for representations associated with classes lying in the same sheet of $G$, providing two alternative formulations. We underline some evidence and illustrate potential consequences.
\end{abstract}

\section{Introduction}


The representation theory of the quantized enveloping algebra $\Ueg$ of a Lie algebra $\g$ at the roots of unity is not completely understood. Big steps towards its comprehension have been made by De Concini, Kac and Procesi in the early 90'es. They have shown that simple modules are always finite-dimensional, they settled a relation between irreducible representations and conjugacy classes in a suitable group $G$ whose Lie algebra is $\g$ by passing through central characters, and they formulated a conjecture relating the dimensions of irreducible representations on the one hand  and of the associated $G$-conjugacy classes on the other hand. The conjecture has been proved in several cases, such as regular classes (\cite{DCKP-regular}), subregular unipotent classes in $SL_n({\mathbb C})$ (\cite{subregular}),  when the order of the root of unity is an odd prime (\cite{modp}),  spherical conjugacy classes (\cite{spherical}, \cite{ccc}). There is also
a not yet published proof by Kremnitzer valid for all unipotent conjugacy classes. However,  a general proof is not available. 
More recently, the representation theory of these algebras has been studied in  \cite{DCM,DCPRR}, where the analysis of branching rules for representations associated with regular conjugacy classes has been carried over.

Motivated by the above mentioned conjecture among other reasons, an analysis of sheets of conjugacy classes has been started in \cite{gio-espo}. These are the irreducible components of the locally closed subsets of $G$ consisting of elements whose conjugacy  class has a fixed dimension.  The main source of inspiration was the analogous work, for adjoint orbits, of Borho and Kraft (\cite{bo,BK}). It is natural to expect that representations of the quantized universal enveloping algebra associated with classes lying in the same sheet should share some properties. 

We recall that $G$ can be parted into a finite union of so-called Jordan classes or decomposition classes. These are irreducible locally closed sets given by  unions  of conjugacy classes with same unipotent part and semisimple part with same connected centralizer. Every sheet contains a dense Jordan class. As a consequence of the reduction theorem in \cite{DCK}, we establish a relation between irreducible representations associated to conjugacy classes in the same Jordan class. It states that there is a natural dimension preserving bijection between the sets of representations associated to conjugacy classes lying in the same Jordan class (Theorem~\ref{thm:jordan}). 

A key role in the description of sheets is played by the induction procedure that produces a conjugacy class in $G$ starting from a conjugacy class in a Levi subgroup of a parabolic subgroup of $G$.  For unipotent conjugacy classes this construction was introduced by Lusztig and Spaltenstein in \cite{lusp}. The dimension of the induced class is related by a simple formula to the dimension of the class in the Levi subgroup. 

Sheets are described as a union of induced conjugacy classes from a Levi subgroup $L$ strictly related to the dense Jordan class. 

On the other hand, a quantum analogue of parabolic subalgebras and of their Levi subalgebras $\l$ are available, and parabolic induction allows for the construction of a $\Ueg$-module from a $\Uel$-module, with good control on central characters. Again, the dimension of the induced module is related by a simple formula to the dimension of the module one is inducing from. For the classes in the sheet lying in the dense Jordan class, induction coincides with saturation and we know that in this case the irreducible representations are induced from irreducible representations of a quantized enveloping algebra of a Levi subalgebra of a parabolic subalgebra by \cite{DCK}.  This is precisely the Lie algebra of the subgroup $L$. However, for the classes in the sheet lying in the boundary of the Jordan class, induction can be different from saturation. 
We conjecture that also in this case the two induction processes should match. We see here an analogy with the results in \cite{lusp}, where it is shown that, over the complex numbers, induction of unipotent conjugacy classes behaves well with respect to Springer's correspondence relating unipotent conjugacy classes of a semisimple algebraic group $G$ to the irreducible representations of its Weyl group. We provide two alternative formulations of our conjecture (Theorem~\ref{thm:equivalent}). If confirmed,  it would shed  new light on the representation theory of $\Ueg$ and its validity could also be applied to the study of the De Concini, Kac and Procesi conjecture (Remark~\ref{consequences}). We provide evidence of our conjecture for cases in which the De Concini, Kac and Procesi's one has been verified, such as regular classes and subregular classes in type $A_n$. 

This paper is an expanded version of a lecture given at the conference ``Hopf algebras and tensor categories'', Almer\'ia, July 2011. One of the aims of the talk was to  underline the interplay between existing results on representations of quantum groups  at the roots of unity and new results on  algebraic groups.  Ideas that might be known to some experts can be expressed in the framework of induced conjugacy classes, decomposition classes and sheets of conjugacy classes, whose systematic description has been only recently made available. We think that it is of interest to put such results on representation theory into this developing  context. 

\section{Notation}\label{preliminaries}

Unless otherwise stated, $H$ is a connected reductive, complex algebraic group, and $G$ is a complex semisimple algebraic group whose simple factors are simply-connected. Both for $G$ and $H$,   $T$ is a fixed maximal torus contained in a Borel subgroup $B$, $\Phi$ is the root system associated with $T$, $\Phi^+$ is the set of positive roots relative to $B$, $\Pi=\{\alpha_1,\ldots,\alpha_n\}$ is the set of simple roots. 
The Borel subgroup opposite to $B$ will be denoted by $B^-$ and the corresponding unipotent  radicals  will be denoted by $U$ and $U^-$.
We shall denote by $Q={\mathbb Z}\Phi$ the root lattice,  by $\Lambda$ the weight lattice, and   
by ${\mathcal W}$ the Weyl group of $G$ or $H$. A standard parabolic subgroup of $G$ or $H$  will be usually denoted by $P$. By $L$ we will indicate a standard Levi subgroup of it, whereas its unipotent radical will be indicated by $U_P$.  By abuse of terminology, we will call $L$ a Levi subgroup of $G$ or of $H$.
By $P^-$ and $U_P^-$ we shall indicate the standard parabolic subgroup opposite to $P$ and its unipotent radical $U^-_{P}$, respectively. The subset of $\Pi$ associated with $L$ will be usually denoted by $\Pi_L$ and we put $B_L=B\cap L$, $B_L^-=B^-\cap L$. The conjugacy class of an element $k$ in a group $K$ will be usually denoted by $\mO_k^K$ and we will denote the operation of conjugation of $x$ on $k$ by $x\cdot k=xkx^{-1}$. The centralizer of $k$ in $K$ will be denoted by $K^k$, whereas the centralizer in $K$ of a subgroup $S$ will be denoted by $C_K(S)$. For any algebraic group $K$, its identity component will be indicated by $K^\circ$.

By a gothic letter we will usually indicate the Lie algebra of the group indicated by the corresponding capital letter, for instance $\g={\rm Lie}(G)$, $\h={\rm Lie}(H)$. 
For an associative algebra $A$, let ${\rm Rep}(A)$ and $\spec(A)$ denote the set of isomorphism classes of  $A$-modules,  and of simple of $A$-modules, respectively. 

\medskip

Let $\ell$ be a positive odd integer coprime with the bad primes of $\g={\rm Lie}(G)$ (cf. \cite[E-12, \S 4]{131}), let $\varepsilon$ be a primitive $\ell$-th root of $1$ and let  $\Ueg$ be the De Concini-Kac specialization (\cite{DCK1}) of the quantized enveloping algebra of $\g$ at $\varepsilon$ corresponding to the isogeny class of $G$. More precisely, the Cartan part $\Ueh$ of $\Ueg$ will be generated by the elements $K_{\xi_i}$ for $i=1,\,\ldots,\,n$ with $\{\xi_i\}_{1\leq i\leq n}$ a basis of  the character group $M=X(T)$  of $T$. 

As usual, we construct the root vectors in $\Ueg$ starting from a fixed reduced decomposition of the longest element $w_0=s_{i_1}\cdots s_{i_N}$ in ${\mathcal W}$ and the associated ordering of the positive roots:
$$\beta_{i_r}=s_{i_1}\cdots s_{i_{r-1}}(\alpha_{i_r}),\quad\mbox{ for } r=1,\ldots,\, N.$$
We consider Lusztig's action of the braid group and the automorphisms $T_i$ (\cite{lu-braid}) so that 
the positive root vectors are then defined as $E_{\beta_r}=T_{i_1}\cdots T_{i_{r-1}}(E_{i_r})$ and the negative ones by $F_{\beta_r}=\omega(E_{\beta_r})$, where $\omega$ is the anti-automorphism interchanging $E_i$ with $F_i$. Then, the quantum version of Poincar\'e-Birkhoff-Witt theorem states that 
$$\{F_{\beta_N}^{a_N}\cdots F_{\beta_1}^{a_1}K_{\xi_1}^{c_1}\cdots K_{\xi_n}^{c_n} E_{\beta_1}^{b_1}\cdots E_{\beta_N}^{b_N},~|~a_i,\,b_j\in{\mathbb Z}_{\geq0};\; c_i\in{\mathbb Z}\}$$ is a basis of $\Ueg$.

It is well-known that the Hopf algebra $\Ueg$ has a large center, containing the Hopf subalgebra $Z_0(G)$ generated by the $\ell$-th powers of the root vectors and of the Cartan generators $K_\lambda$. This subalgebra is independent of the choice of the reduced expression of the longest element $w_0$ of ${\mathcal W}$, as it is the minimal subalgebra of $\Ueg$ closed under the Poisson bracket and containing the $\ell$-th powers of the Chevalley generators and of the $K_\lambda$, for $\lambda\in M$ (\cite[Page 22]{DCM}).

Restriction of an irreducible $\Ueg$-module $V$ to $Z_0(G)$ determines a natural map $\Chi\colon \spec(\Ueg)\to\spec(Z_0(G))$ obtained by associating to $V$ its $Z_0(G)$-central character.  


In \cite[\S 4]{DCKP} a natural map $\pi\colon \spec(Z_0(G))\to B^-B\subset G$ is defined and it is shown that $\pi$ is an unramified covering of the big cell $B^-B$. 
It is obtained as follows.
Since $Z_0(G)\cong Z_0^-(G)\otimes Z_0^0(G)\otimes Z_0^+(G)$ where $Z_0^{\pm}(G)$ are generated by $\ell$-th powers of root vectors of fixed positivity and $Z_0^0(G)$ is generated by the $\ell$-th powers of the Cartan elements $K_\lambda$, we have a decomposition $$\spec(Z_0(G))\cong\spec(Z_0^+(G))\times\spec(Z_0^0(G))\times\spec(Z_0^-(G)).$$ The map $\pi$ is a product of the map $\pi^0\colon \spec(Z_0^0(G))\cong T\to T$ obtained by taking the square of an element, with two birational isomorphisms $$\pi^{\pm}\colon \spec(Z_0^{\pm}(G))\to U^{\pm}$$ which are constructed as follows. Let  $f_i,\,e_i\in\g$ be  Chevalley generators, let $f_\beta\in\g$ be the root vectors constructed by Tits by using the operators $T'_i=\exp({\rm ad} f_i)\exp({\rm ad}e_i)\exp({\rm ad} f_i)$. Let $T_0=T_{i_1}\cdots T_{i_N}$,  and let $T'_0=T'_{i_1}\cdots T'_{i_N}$. 
Then
$$\pi^-(\chi)=\exp(\chi(y_{\beta_N})f_{\beta_N})\cdots \exp(\chi(y_{\beta_1})f_{\beta_1})$$ where $y_{\beta_r}=c_{\beta_r}^\ell F_{\beta_r}^\ell$ for  suitable scalars $c_{\beta_r}$ and
$$\pi^+(\chi)=\exp(\chi(T_0(y_{\beta_N}))T'_0(f_{\beta_N}))\cdots \exp(\chi(T_0(y_{\beta_1}))T'_0(f_{\beta_1})).$$
Then, $\pi(\chi^-,\chi_0,\chi^+)=\pi^-(\chi^-)\pi_0(\chi_0)\pi^+(\chi^+)$. Let us observe that $\chi^+=0$ if and only if  $\pi(\chi)\in B^-$.

It is shown in \cite[\S 7]{DCKP} that, for $G$ simply-connected, the group $\spec(Z_0(G))$ is isomorphic to the Poisson dual group $H_G$ of $G$, where 
$$H_G=\{(sv,tu)\in TU^-\times TU~|~t=s^{-1}\}.$$
After this identification, the map $\pi$ is given by $(sv,tu)\mapsto (vs)^{-1}(tu)\in B^-B$. 

\begin{remark}\label{rk:properties}
The composition  $\pi\circ \Chi$ has the following properties:

\begin{enumerate}
\item The fiber of $g\in B^-B$ through $\pi\circ \Chi$ is $\spec(U_\chi(G))$ for some fixed finite-dimensional quotient of $\Ueg$, namely ${\cal U}_\chi(G)\cong\Ueg/(z-\chi(z),z\in Z_0(G))$, with $\pi(\chi)=g$  (\cite[Thm 6.1]{DCKP}). 
\item The Poincar\'e-Birkhoff-Witt theorem is compatible with taking the quotient by $(z-\chi(z),z\in Z_0(G))$ so
$$\{F_{\beta_N}^{a_N}\cdots F_{\beta_1}^{a_1}K_{\xi_1}^{c_1}\cdots K_{\xi_n}^{c_n} E_{\beta_1}^{b_1}\cdots E_{\beta_N}^{b_N},~|~ 0\leq a_i,\,b_j,\,c_i\leq \ell-1\}$$ is a basis of $\Uegchi$.
\item If $g=\pi(\chi)$ and $h=\pi(\xi)$ are $G$-conjugate and lie in $B^-B$, then $\Uegchi\cong \Uegxi$  (\cite[Thm 6.6]{DCKP}). Therefore, it is not restrictive to look only at representations of $\Uegchi$ with $\pi(\chi)\in B^-$.
\end{enumerate}
\end{remark}

As a consequence of the above remark, a map $\varphi$ from $\spec(\Ueg)$ to the set of conjugacy classes of $G$ is defined. 

De Concini, Kac and Procesi formulated in \cite{DCKP} the following conjecture: 

\begin{conjecture}\label{conj:DCKP}If $V\in\spec(\Ueg)$ and $\mO=\varphi(V)$ then $\ell^{\frac{\dim(\mO)}{2}}$ divides $\dim V$.
\end{conjecture}

Conjecture~\ref{conj:DCKP} has been confirmed  in several cases listed in the Introduction. 
It is natural to try to seek for  relations among fibers through $\varphi$ of conjugacy classes of the same dimension in a given family. We will do so for Jordan classes (decomposition classes) and for sheets.

\subsection{Quantized Levi subalgebras and parabolic subalgebras}

Let us consider a standard parabolic subgroup $P$ of $G$ with standard Levi decomposition $P=LU_P$ and basis $\Pi_L\subset \Pi$ of the corresponding root subsystem $\Phi_L$. We recall that $L$ is simply-connected since $G$ is so.

The subalgebra generated by $\Ueh$ and the $E_i$, $F_i$ for $\alpha_i\in\Pi_L$ is isomorphic to $\Uel$.  

Let $w_L$ be the longest element in the Weyl group ${\mathcal W}_L$ of $L$ and let $N_L=|\Phi_L\cap\Phi^+|$. From now on we shall always consider a reduced decomposition of $w_0=s_{i_1}\cdots s_{i_N}$ such that the product of the first $N_L$ terms is equal to $w_\Pi$. Thus, $E_\beta$ and $F_\beta$ lie in $\Uel$ for every $\beta\in\Phi_L^+$ and the PBW-bases for $\Uel$ and $\Ueg$ are compatible. Like we did for $G$, we may define $Z_0(L)\subset Z_0(G)$ and the restriction map $\Chi_L$. Moreover,  
the construction of $\pi$ and the dual group holds for any connected reductive algebraic group with simply-connected simple factors. In particular, it holds for $L$ a Levi subgroup of $G$ (\cite[\S 5]{DCM}) and it is clear from its description that the map $\pi$ and the corresponding map $\pi_L$ on the Poisson dual $H_L$ are compatible. Besides, by our choice of the reduced decomposition of $w_0$, if $\chi\in \spec (Z_0(G))$ and $\chi^+=0$, then for the corresponding maps $\pi_L^{\pm}$, the restriction to $Z_0(L)$ of $\chi^+$ is $\chi^+_L=0$ and $\pi^-(\chi^-)\in U^-_{P}\pi_L^-(\chi^-_L)$.

An analogue of the map $\varphi$ can be defined also on $\spec(\Uel)$ and we will indicate it by $\varphi_L$.  

%

By \cite[Lemma 5.1]{DCM}  we have 
\begin{equation}\label{eq:DCM}\Uel\cong ({\cal U}_\varepsilon([L,L])\otimes k[Z(L)^\circ])^\Gamma\end{equation}
where $\Gamma$ is the kernel of the isogeny $\phi\colon [L,L]\times Z(L)^\circ \to L$ and its action on ${\cal U}_\varepsilon([L,L])\otimes k[Z(L)^\circ]$ is the natural one.

%
The quantized parabolic subalgebra $\Uep$ is the subalgebra of $\Ueg$ generated by $\Uel$ and all the $E_i$'s.
By the Poincar\'e-Birkhoff-Witt theorem, $\Uep$ is again finite over the subalgebra $Z_0(P)$ generated by $Z_0(L)$ and all the $\ell$-th powers of the positive root vectors $E_\alpha$. Let us consider $W\in \spec(\Uep)$, such that its associated $Z_0(P)$-character $\chi$ is trivial on $Z_0(P)\cap Z_0^+(G)$. Then, the ideal of $\Uep$ generated by the root vectors $E_\alpha$ for $\alpha$ in $\Phi^+\setminus\Phi_L$ is nilpotent, so it must act trivially on $W$. Hence, $W\in \spec(\Uel)$, with $\Uel$ viewed as a quotient of $\Uep$. Conversely, the action on every $\Uel$-module may be extended uniquely to an irreducible $\Uep$-action by letting the $E_\alpha$ for $\alpha$ in $\Phi^+\setminus\Phi_L$ act trivially. In other words, if $\chi\in\spec(Z_0(L))$ is such that $\chi^+=0$, then we may consider the quotients
$\Uelchi=\Uel/(z-\chi(z),z\in Z_0(L))$ and
$\Uepchi=\Uep/(E_\alpha^\ell, \alpha\in\Phi^+;\; z-\chi(z),z\in Z_0(L))$ and
there is a natural  bijection between $\spec(\Uelchi)$ and $\spec(\Uepchi)$.

Let $\chi$ be a central character in $\spec(Z_0(G))$ with $\pi(\chi)\in B^-$ and let $\chi_L$ be its restriction to $\zol$. We define the parabolic induction map 
$$\begin{array}{rl}
{\mathcal Ind}_L^{G,\chi}\colon \spec(\Uelchil)&\to{\rm Rep}(\Uegchi)\\ 
V&\to \Uegchi\otimes_{\Uepchil} V\end{array}$$ 

By Poincar\'e-Birkhoff-Witt's theorem 
we have 
\begin{equation}\label{eq:dime}\dim  {\mathcal Ind}_L^{G,\chi}V=\ell^{|\Phi^+|-|\Phi_L^+|}\dim V.\end{equation}

Conversely, if $V\in\spec(\Uegchi)$ with $\pi(\chi)=b^-\in B^-$, and if $V={\mathcal Ind}_{L}^{G,\chi}W$ for some $W\in \spec(\Uep)$ and some standard parabolic subgroup $P$ with Levi subgroup $L$, then  $\Chi_L(W)$ is obtained by restriction of $\Chi(V)$ to $Z_0(L)$.  By construction, $\pi\circ \Chi(V)\in (U^-_{P})(\pi_L\circ X_L)(W)$ so we have $\pi_L\circ X_L(W)=b'\in B\cap L$ and  $b\in U^-_{P}b'$ .

\subsection{Jordan classes}\label{sec:jordan}

We are ready to describe the relation with the first type of families of conjugacy classes. Notation is as in Section~\ref{preliminaries}.

\begin{definition}Let $H$ be a connected reductive algebraic group. A {\em Jordan class} or {\em decomposition class} is an equivalence class with respect to the relation on $H$ defined by: $g\sim h$ if 
$g=su$, $h=rv$ and there exists $x\in H$ such that: $H^{xsx^{-1}\circ}=H^{r\circ}$; the classes $\mO_{xux^{-1}}^{H^{r\circ}}$ and $\mO_v^{H^{r\circ}}$ coincide; and  
$xsx^{-1}\in Z(H^{r\circ})^\circ r$. 
\end{definition}

There are finitely many Jordan classes in $H$ and the Jordan class $J(su)$  of $g=su$ equals $H\cdot(Z(H^{s\circ})^\circ s)^{reg}u$ where by  $(Z(H^{s\circ})^\circ s)^{reg}$ we denote the subset of $Z(H^{s\circ})^\circ s$ consisting of elements whose centralizer has minimal dimension, that is, of those elements in $Z(H^{s\circ})^\circ s$ whose connected centralizer equals $H^{s\circ}$. 

In general, $H^{s\circ}$ is not necessarily a Levi subgroup of $H$ but $L=C_H(Z(H^{s\circ})^\circ)$ is so. It is the minimal Levi subgroup of $H$ containing $H^{s\circ}$ and we call it the {\em Levi envelope} of $H^{s\circ}$. There holds $Z(L)^\circ= Z(H^{s\circ})^\circ$.

Let us also recall that if $H$ is simply-connected, then the centralizer of any semi-simple element is connected.

We reformulate the following well-known result:

\begin{theorem}\label{thm:DCK-reduction}(\cite[Theorem \S 8]{DCK}) Let $G$ be semi-simple and simply connected. Let $g=us\in U^-T$ be the Jordan decomposition of an element in $G$ such that $L=C_G(Z(G^{s})^\circ)$ is a {\em proper} standard Levi subgroup of $G$. Then, for every $V$ in $\varphi^{-1}(\mO_g^G)$, with $\pi(\chi)=g$,  there exists a unique $W\in\varphi_L^{-1}(\mO_g^L)$ such that $V={\mathcal Ind}_L^{G,\chi}W$. Besides, parabolic induction establishes a bijection between $\spec(\Uegchi)$ and $\spec(\Uelchi)$. 
\end{theorem} 

As a consequence of the above Theorem, we have:

\begin{theorem}\label{thm:jordan}Let the conjugacy classes $\mO_g^G$ and $\mO_h^G$ lie in  the same Jordan class, with $g$ having Jordan decomposition $g=us\in U^-T$ and  such that
the Levi subgroup $L=C_G(Z(G^{s})^\circ)$ is standard. Let $\pi(\chi)=g$ and $\pi(\xi)=h$. Let $\varphi^{-1}(\mO_g^G)=\{V_1,\ldots,\,V_m\}$ and 
$\varphi^{-1}(\mO_h^G)=\{V'_1,\ldots,\,V'_n\}$. Then, $m=n$ and there is a $1$-dimensional $\Uel$-module $V_\lambda$ such that, up to a permutation of the indices, if $V_i={\mathcal Ind}_L^{G,\chi}W_i$ then $V'_i={\mathcal Ind}_L^{G,\xi}(W_i\otimes V_\lambda)$. \end{theorem}
\pf Without loss of generality we may assume that the Jordan decomposition of $h$ is $h=ur$ with $G^{r}=G^{s}$ and 
$r=zs$ for $z\in Z(L)^\circ$. A central character $\chi_z$ such that $\pi(\chi_z)=z$  satisfies $\chi_z(K_\mu^{2\ell})=\mu(z)$ for every $\mu\in \Lambda$. In particular, for every $\alpha\in\Phi_L$ we have $\chi_z(K_\alpha^{2\ell})=1$. We may assume thus that $\xi=\chi\chi_z$.

By \eqref{eq:DCM}, the subalgebra $\Uel$ is a subalgebra of the tensor product ${\cal U}_\vep([L,L])\otimes {\mathbb C}[Z(L)^\circ]$.
Moreover, $Z(L)^\circ= M_L\otimes_{\mathbb Z}{\mathbb C^*}$ where $M_L$ is the lattice of cocharacters which are trivial on $\Pi_L$. 
We shall denote by $\{\eta_1^\vee,\,\ldots,\eta_k^\vee\}$ a basis of $M_L$, for $k=|\Pi|-|\Pi_L|$, and by  $\{\theta_1,\,\ldots,\,\theta_k\}$ its dual basis in ${\rm Hom}_{\mathbb Z}(M_L,{\mathbb Z})$.

Let us fix scalars $\lambda_i$ for $i=1,\,\ldots,\,k$  satisfying $\lambda_i^{2\ell}=\theta_i(z)$. Such scalars define a one-dimensional representation of ${\mathbb C}[Z(L)^\circ]\cong {\mathbb C}[K^{\pm1}_{\theta_i}]_{1\leq i\leq k}$ which we shall denote by  ${\mathbb C}_\lambda$. Tensoring with the trivial ${\cal U}_\vep([L,L])$-module, we obtain a ${\cal U}_\vep([L,L])\otimes {\mathbb C}[Z(L)^\circ]$-module $V_\lambda$ whose restriction to $\Uel$ is a one-dimensional module with central character  $\chi_z$.

Let $V_i$ be an irreducible $\Uegchi$-module. Then, there is a unique irreducible $\Uelchi$-module $W_i$ such that 
$V_i={\mathcal Ind}_L^{G,\chi}W_i$. On the other hand, tensoring with $V_\lambda$ sets up a bijection between $\spec(\Uelchi)$ and $\spec({\cal U}_{\chi\chi_z}(L))$. Parabolic induction determines a bijection between $\spec({\cal U}_{\chi\chi_z}(G))$ and $\spec({\cal U}_{\chi\chi_z}(L))$. Hence, $m=n$ and, up to reordering of the terms, we have $V_i'={\mathcal Ind}_{L}^{G, \chi\chi_z}(W_i\otimes V_\lambda)$ for every $i=1,\,\ldots,\,m$. \hfill$\Box$

\medskip

It follows from the above theorem, although it is already implicit in Theorem~\ref{thm:DCK-reduction}, that in order to verify Conjecture~\ref{conj:DCKP} it is enough to confirm it for the finitely many Jordan classes, because of \eqref{eq:dime}. 

It is pointed out in \cite[Remark 8.1]{DCK} that the $\Uel$-module $W$ in  Theorem~\ref{thm:DCK-reduction} remains irreducible when restricted to the subalgebra $\Uel'$ generated by the root vectors $E_\alpha,\,F_\beta$ corresponding to roots in $\Pi_L$ and by the $K_\beta$ for $\beta$ in the root lattice of $\Phi_L$. with notation as in Theorem~\ref{thm:jordan}, the irreducible $\Uel'$-modules associated with $V_i$ and $V_i'$ coincide.

\subsection{Sheets of conjugacy classes}\label{sec:sheets}

In this section we deal with the second type of families of conjugacy classes.

Let $H$ be a connected reductive algebraic group and let $X$ be a $H$-variety. The irreducible components of the locally closed subsets
$$X_{(n)}=\{x\in X~|~\dim H\cdot x=n\}$$
are called the sheets of $X$. The sheets for the adjoint action on  $X=\h$ have been studied in detail in  \cite{bo,BK}. In analogy to this situation, the sheets of conjugacy classes in $X=H$ have been studied in \cite{gio-espo}. Every sheet in $H_{(n)}$ can be described as $H_{(n)}\cap \overline{J(g)}$ for a unique Jordan class $J(g)$.  A sheet $S$ whose dense Jordan class is $J(g)$ will be also denoted by $S(g)$. Clearly, if $g\sim h$, then $S(g)=S(h)$. 

When we write $S=S(su)$ we shall always mean that $su$ is the Jordan decomposition of an element in $H$ and that $s,\,u$ are chosen with $u\in U^-$ and $s\in T$ such that the Levi-envelope of $H^{s\circ}$ is standard. Sheets are best described in terms of  {\em induced conjugacy classes}. These are defined as follows.

Let  $L$ be a Levi subgroup in $H$, let $P=LU_P$ be a parabolic subgroup of $H$, and let $\mO_l^L$ be a conjugacy class in $L$. The $H$-{\em conjugacy class induced by $\mO_l^L$} is $\mO:={\rm Ind}_{L,P}^H(\mO_l^L):=H\cdot(\mO_l^LU_P)^{reg}$, where $reg$  denotes the subset of elements with centralizer of minimal dimension. Since the semisimple parts of the elements in  $\mO_l^L U_P$ are all $H$-conjugate and, by work of Dynkin and \cite[Corollary 3.7, Lemma 5.1]{kostant}, there are only finitely many unipotent classes in a reductive group, the set $\mO$ is indeed a $H$-conjugacy class. 
Induced unipotent conjugacy classes have been extensively studied in \cite{lusp}, whereas induced adjoint orbits have been addressed in \cite{bo}.  In particular, they are independent  of the choice of the parabolic subgroup $P$ with  Levi factor $L$.  It was shown in \cite{gio-espo} that 
\begin{equation}\label{indotta}{\rm Ind}_{L,P}^H(\mO_{su}^L)=H\cdot(s\,{\rm Ind}_{L^{s\circ}}^{H^{s\circ}}(\mO_u^{L^{s\circ}})).\end{equation}
Hence, independence of the parabolic subgroup, transitivity, and the dimension formula in \cite{lusp} follow in the general group case from the unipotent case:
\begin{equation}\label{dimensione}\dim {\rm Ind}_{L,P}^H(\mO_{su}^L)=\dim H-\dim L+\dim \mO_{su}^L\end{equation} 
and we may omit the index $P$  in ${\rm Ind}^H_{L,P}$.

For a sheet $S=S(su)$ we have: 
\begin{equation}\label{sheet}S=H_{(n)}\cap\overline{J(su)}=\bigcup_{z\in Z(H^{s\circ})^\circ}{\rm Ind}_L^H(\mO_{zsu}^L)\end{equation}
where $L$ is the Levi-envelope of $H^{s\circ}$. 

The conjugacy classes in $S$ lying in the dense subset $J(su)$ are those for which $zs\in (Z(H^{s\circ})^\circ s)^{reg}$ and in this case 
 ${\rm Ind}_{L}^H(\mO_{ru}^L)=\mO^H_{ru}$.
For those classes,  Theorem~\ref{thm:jordan} establishes a relation between fibers of $\varphi$. We would like to analyze the situation at boundary points, that is, at classes in $(\overline{J(su)}\setminus J(su))\cap S$. Let us first state a few basic properties about parabolic induction of representations and
induction of conjugacy classes. 

\begin{lemma}\label{transitivity}Let $G$ be simply-connected. Let $V\in\spec(\Ueg)$, let $L\subset L'$ be standard Levi subgroups of the standard parabolic subgroups $P\subset P'$. Let $\chi\in\spec(\zog)$ with $\pi(\chi)\in B^-$ and  assume $V={\mathcal Ind}_{L}^{G,\chi}W$ for some $W\in \spec(\Uelchil)$. Then $V={\mathcal Ind}_{L'}^{G,\chi}W'$ for some $W'\in \spec({\cal U}_{\chi_{L'}}(L'))$. If, in addition, $\mO_{\pi(\chi)}^G={\rm Ind}_{L}^G(\mO_{\pi_L(\chi_L)}^L)$
then 
$$\mO_{\pi(\chi)}^G={\rm Ind}_{L'}^G(\mO_{\pi_{L'}(\chi_{L'})}^{L'}),\quad\mbox{and}\quad \mO_{\pi(\chi_{L'})}^{L'}={\rm Ind}_{L}^{L'}(\mO_{\pi_L(\chi_L)}^L).$$
\end{lemma}
\pf The first assertion is standard:
let $X=\Uegchi$, $Y={\cal U}_{\chi_{L'}}(P')$, $Z=\Uepchil$.
Then $V\cong X\otimes_Z W\cong X\otimes_Y Y\otimes_Z W$ and the $Z$-module $W'=Z\otimes_Y W$ is necessarily irreducible because $V$ is so. 

For the second statement we have $\pi(\chi)=ul\in(U^-_{P}l)^{reg}$, and $u=vl'\in U^-_{P}\cap U^-_{P'}L'$, with $l=\pi(\chi_L)$ and $l'l=\pi(\chi_{L'})$. Therefore we need to show that $ul=vl'l\in (U^-_{P'}l'l)^{reg}$ in $G$ and 
$l'l\in (U^-_{P})^{reg}$ in $L'$. It is enough to prove that
$\dim \mO_{ul}^G=\dim {\rm Ind}_{L'}^G(\mO_{l'l}^{L'})$ 
and $\dim \mO_{l'l}^{L'}=\dim{\rm Ind}_L^{L'}(\mO_l^L)$. We have
%
%
$$\begin{array}{rl}
\dim G-\dim L+\dim \mO_l^L &=\dim \mO_{ul}^G\\
&=\dim \mO^G_{vl'l}\\
&\leq\dim{\rm Ind}_{L'}^G(\mO_{l'l}^{L'})\\
&= \dim \mO^{L'}_{l'l}+\dim G-\dim L'\\
&\leq \dim{\rm Ind}_{L}^{L'}(\mO_l^L)+\dim G-\dim L'\\
&=\dim L'-\dim L+\dim \mO_l^L +\dim G-\dim L'\\
&=\dim G-\dim L+\dim \mO_l^L 
\end{array}$$ so equality holds in all steps, yielding the statement. 
\hfill$\Box$

\begin{proposition}\label{saturated}Let $L$ be a standard Levi subgroup of a connected reductive algebraic group $H$ and let 
$su$ be the Jordan decomposition of a representative of the $L$-class $\mO_{su}^L$, with 
$s\in T$ and $u\in U\cap L^{s\circ}$. Then, ${\rm Ind}_L^H(\mO^L_{su})=\mO_{su}^H$ if and only if $L\supset H^{s\circ}$.
\end{proposition}
\pf If $L\supset H^{s\circ}$, equation \eqref{indotta} shows that ${\rm Ind}_L^H(\mO_{su}^L)=\mO_{su}^H$. Conversely, if  equality holds, 
\eqref{dimensione} and \eqref{indotta}  give $\dim H^{su}=\dim L^{su}$ so $H^{su\circ}=L^{su\circ}$. 

We shall show that this is possible only if $H^{s\circ}\subset L^{s\circ}$. Let   
 $\Phi_s=\{\alpha\in\Phi~|~\alpha(s)=1\}$. We may thus write:
$$H^{s\circ}=\langle T, X_\alpha,\;\alpha\in\Phi_s\rangle;\quad L=\langle T,\,X_\alpha,\;\alpha\in \Phi_L\rangle.$$

We may choose a basis for $\Phi_s$ consisting of not necessarily simple, but positive roots. Let $\beta$ be a highest root in a simple component  $\Psi$ of $\Phi_s$. 
Then, $X_\beta$ commutes with $u\in H^{s\circ}\cap U$, so it lies in $H^{su,\circ}=L^{su,\circ}$. Since $L$ is standard, the support of $\beta$ lies in $\Phi_L$. Thus, the basis of each component of  $\Phi_s$ lies in $\Phi_L$. \hfill$\Box$

\begin{proposition}Let $G$ be simply-connected, let $s\in T$ and let $L$ be a standard Levi subgroup such that $G^{s}\subset L$. Let $V$ be an irreducible $\Uegchi$-module such that $\varphi(V)=\mO_{su}^G$. Then $V={\mathcal Ind}_L^{G,\chi}W$ for a unique $W$ such that $\varphi_L(W)=\mO_{su}^L$.
\end{proposition}
\pf By Lemma~\ref{transitivity} and Proposition \ref{saturated}, in order to prove the existence of $W$, it is enough to consider the case of $L$ a minimal standard Levi subgroup containing $G^{s}$, that is,  when $L$ is the Levi-envelope of $G^s$. This is Theorem~\ref{thm:DCK-reduction}. The proof of uniqueness in there applies also to the case of a general standard Levi subgroup $L$ containing $G^s$. We recall here the argument for completeness. Let $V={\mathcal Ind}_L^{G,\chi}W={\mathcal Ind}_L^{G,\chi}W'$ for two irreducible $\Uelchi$-modules $W$ and $W'$.
By Remark~\ref{rk:properties}, part 3 and Proposition~\ref{saturated}, we may assume that $\pi(\chi)=\pi(\chi_L)$, that is, we may assume that $\chi(F_\alpha^\ell)=0$ for every $\alpha\not\in\Phi_L$. Thus, $\Uegchi$ has a unique natural  ${\mathbb Z}$-grading obtained by setting $\deg(\Uelchil)=0$ and $\deg(F_\alpha)=1$ for $\alpha\in\Pi\setminus\Pi_L$. The $\Uegchi$-module $V$ is naturally ${\mathbb Z}$-graded by setting $V_0=W$. By construction, $V_j=0$ for $j<0$.  It is not hard to prove that the natural projection $\pi_0$ of $V$ onto $V_0=W$  is $\Uelchil$-equivariant. Therefore,  its restriction to $W'$ is either an isomorphism or the trivial map. However, if $\pi_0(W')=0$, that is, if $W'\subseteq\bigoplus_{j>0}V_j$,  then we would have 
$$V={\mathcal Ind}_L^{G,\chi}W'=\sum_{a_i\geq0}F_{\beta_N}^{a_N}\cdots\,F_{\beta_{N_L+1}}^{a_{N_L+1}}W'\subseteq\bigoplus_{j>0}V_j,$$  which is impossible.
\hfill$\Box$

\bigskip

The above statements lead us to the formulation of the following conjecture:

\begin{conjecture}\label{conj:gio}Let $G$ be simply-connected. Let $V$ be an element of $\spec(\Uegchi)$ with $\pi(\chi)=g\in B^-$. 
If $\mO_g^G={\rm Ind}_{L'}^G(\mO_l^{L'})$ for some class $\mO_l^{L'}$ in a standard Levi subgroup $L'$ of a parabolic subgroup $P$ of $G$, then, there exists $W\in\spec({\cal U}_{\chi_{L'}}(L'))$  such that $V={\mathcal Ind}_{L'}^{G,\chi}(W)$.  
\end{conjecture}

We point out that, by construction, if a module $W$ such that $V={\mathcal Ind}_{L'}^{G,\chi}(W)$ exists,  then $W$ is irreducible and we necessarily have $\varphi_{L'}(W)=\mO_{l}^{L'}$. 

\bigskip

Let $S=S(su)$ and let $\mO^G_h\in J(su)$. Then, by \eqref{indotta} we have the equality $\mO_h^G={\rm Ind}_{L}^G(\mO_{su}^L)$ for $L$ the Levi-envelope of $G^{s}$. So, Theorem~\ref{thm:DCK-reduction} confirms Conjecture~\ref{conj:gio} in this case. 
For this reason, Conjecture~\ref{conj:gio} should be seen as an extension to $\overline{J(su)}\cap S$ of Theorem~\ref{thm:DCK-reduction}.  
We will now see that Conjecture~\ref{conj:gio} may be formulated as a statement about sheets.

We recall that an element $h\in H$ with Jordan decomposition $h=su$ is called  {\em isolated} (\cite{lusztig}) or {\em exceptional} (\cite{DCK}, where a complete list of exceptional semisimple elements is given) if the Levi envelope of $H^{s\circ}$ is $H$.

A unipotent element in  $H$ is called {\em rigid} if its class is not induced from a class of any proper Levi subgroup of $H$. Rigid unipotent conjugacy classes have been classified in \cite{ela, ke}, and a complete list is available in \cite{spaltenstein}. Each rigid unipotent class is itself a sheet in $[H,H]$. In analogy to this, rigid elements in a semisimple group $H'$ have been defined in \cite{gio-espo} as those whose class is a single sheet in $H'$. We may extend this definition to elements in $H$ by saying that $h\in H$ is rigid if $Z(H)^\circ\mO_h^H$ is a sheet.  

The following statement can be deduced from \cite{gio-espo} and we state it here for completeness' sake.

\begin{lemma}\label{lem:equivalent} Let $H$ be a connected reductive group and let $h\in H$ with Jordan decomposition $h=su$.
Then the following are equivalent:
\begin{enumerate}
\item $\mO_h^H$ is rigid;
\item The element $h$ is isolated in $H$ and $\mO^{H^{s\circ}}_u$ is a rigid unipotent class in $H^{s\circ}$;
\item $\mO_{h}^H$ is not induced from a class in any proper Levi subgroup of $H$.
\end{enumerate}
\end{lemma}
\pf Assume $1$ holds.  Then, $J(h)\subset Z(H)^\circ \mO_h^H\subset J(h)$ so, if $h=su$ is the Jordan decomposition of $h$, we have 
$Z(H^{s\circ})^\circ=Z(H)^\circ$ that is, $h$ is isolated. Moreover, by \cite[Thm. 5.6(a)]{gio-espo}, the unipotent class $\mO_u^{H^{s\circ}}$ is rigid, whence $2$ holds. 

Assume now that $2$ holds. Then, by \cite[Prop. 5.3(b,c)]{gio-espo}, the class $\mO_{h}^H$ is not induced from a class in any proper Levi subgroup of $H$.

Assume now $3$ holds and let $S=S(rv)$  be a sheet containing $\mO_{h}^H$. Then by \eqref{sheet} we necessarily have $C_H(Z(H^{r\circ}))=H$, and, up to conjugation, $su=zrv$ for some $z\in Z(H)^\circ$ so $S=Z(H)^\circ\mO^H_{rv}=Z(H)^\circ\mO_{su}^H$.\hfill$\Box$

\begin{theorem}\label{thm:equivalent}Let $G$ be simply-connected, let $V$ be an element of $\spec(\Uegchi)$ and let $\pi(\chi)=g\in B^-$. The following are equivalent:
\begin{enumerate}
\item Conjecture~\ref{conj:gio}.
\item If $\mO_g$ lies in the sheet $S=S(us)$, then, there exists $W\in\spec(\Uelchil)$  such that $V={\mathcal Ind}_L^{G,\chi}(W)$, for  $L$  the Levi-envelope of $G^{s}$ and $\pi_L(\chi_L)=us$.  
\end{enumerate}
\end{theorem}
\pf By \eqref{sheet}, Conjecture~\ref{conj:gio} implies statement $2$. On the other hand, Lemma~\ref{transitivity} shows that it is enough to confirm Conjecture~\ref{conj:gio}  for $\mO_g^G={\rm Ind}_L^G(\mO_l^L)$ and $\mO_l^L$ rigid in $L$. By Lemma~\ref{lem:equivalent}, the element $l=us$ is isolated and $\mO_{u}^{G^{s}}$ is a rigid unipotent class. In other words, the Levi-envelope of $G^{s}$ is $L$ so by
 \cite[Thm. 5.6(a)]{gio-espo} and \eqref{sheet}, the Jordan class $J(l)$ is dense in a sheet containing $\mO_g^G$. \hfill$\Box$

\bigskip

%
%
%
%

The remainder of the paper is devoted  to consequences of this conjecture and  evidence of it.


\begin{remark}\label{consequences}Assume Conjecture~\ref{conj:gio} holds for a conjugacy class $\mO_g^G={\rm Ind}_L^G(\mO_l^L)$. Then, if Conjecture~\ref{conj:DCKP} holds for $\mO_l^L$ then Conjecture~\ref{conj:DCKP} holds for $\mO_g^G$. Indeed, by \eqref{eq:dime} 
$$\dim V=\dim {\mathcal Ind}_L^{G,\chi}W=\ell^{|\Phi^+|-|\Phi_L^+|}\dim W=\ell^{\frac{\dim G-\dim L}{2}} \dim W$$
so if $\ell^{\frac{\dim \mO_l^L}{2}}$ divides $\dim W$ then 
$\ell^{\frac{\dim \mO_l^L+\dim G-\dim L}{2}}=\ell^{\frac{\dim \mO_g^G}{2}}$ divides $\dim V$ by \eqref{dimensione}.\end{remark}

\bigskip

We have already pointed out that in order to prove Conjecture~\ref{conj:DCKP} for a group $G$, it is enough to prove 
it for isolated classes.  If Conjecture~\ref{conj:gio} were confirmed,  in order to prove Conjecture~\ref{conj:DCKP}, it would be enough to 
prove it for all rigid classes.  Since in type $A$ all Levi subgroups are of type $A$, isolated elements are only unipotent, and unipotent rigid elements are trivial (\cite{spaltenstein}), Conjecture~\ref{conj:gio} would imply Conjecture~\ref{conj:DCKP} recovering Kremnitzer's result in type $A$.
%
%
%
We can similarly deal with a few more cases. 

\begin{proposition}If $G$ is of type $G_2$, then Conjecture~\ref{conj:gio} for $G$ implies Conjecture~\ref{conj:DCKP} for $G$.
\end{proposition}
\pf By Remark \ref{consequences} we would need to verify Conjecture~\ref{conj:DCKP} for rigid conjugacy classes in $G$ and rigid conjugacy classes in all possible Levi subgroups of $G$.
The rigid conjugacy classes in $G$ are: the two isolated semisimple classes in the list in \cite{DCK}, and the unipotent classes labeled by $0$, $\tilde{A}_1$ and $A_1$ according to Elashvili's list of rigid classes (\cite[p. 173]{spaltenstein}). They are all spherical so Conjecture~\ref{conj:DCKP} holds for representations associated with such classes by \cite{ccc}. Proper Levi subgroups in $G$ are all of type $A_1$, so the statement follows.\hfill$\Box$

\begin{proposition}If $G$ is of type $C_3$, $B_3$, or $D_4$, then Conjecture~\ref{conj:gio} for $G$ implies Conjecture~\ref{conj:DCKP} for $G$.
\end{proposition}
\pf  In type $C_3$ the rigid orbits are either isolated semisimple, unipotent of type $0$ or $A_1$, nor semisimple nor unipotent with centralizer of type $A_1\times C_2$ and unipotent part in $C_2$ of type $A_1$. Such classes are all spherical, hence Conjecture~\ref{conj:DCKP} holds for them by \cite{ccc}. Since Levi subgroups are either of type $C_2$ or products of type $A$, and since in type $C_2$ Conjecture~\ref{conj:DCKP} is confirmed in \cite[Corollary 33]{ccc}, the statement follows for $C_3$.

In type $B_3$ the non-trivial rigid classes are either semisimple isolated or of type $2A_1$, hence they are all spherical and Conjecture~\ref{conj:DCKP} is confirmed by \cite{ccc}. Proper Levi subgroups are of type $C_2$ or products of type $A$, so if Conjecture~\ref{conj:gio} were confirmed,  Conjecture~\ref{conj:DCKP} would be verified in type $B_3$.

In type $D_4$ non-trivial rigid classes are either rigid unipotent with partition $(3,2^2,1)$ or $(2^2,1^4)$ or semisimple isolated, and they are spherical, so Conjecture~\ref{conj:DCKP} holds for them. Since all Levi subgroups are products of type $A$, we have the statement.\hfill$\Box$ 





We end the paper by pointing out evidence of Conjecture~\ref{conj:gio} available in the literature.
%

\begin{proposition}Let $\mO_g^G$ be a regular conjugacy class in $G$ and let $V\in\spec(\Uegchi)$ with $\pi(V)=\mO_g^G$. Then Conjecture~\ref{conj:gio} holds for $\mO_g^G$.
\end{proposition}
\pf Every regular conjugacy class is induced from the class of its semisimple part, with trivial Levi subgroup  $L=T$. Let $v$ be a highest weight vector in $V$. Then $W={\mathbb C}v$ is an irreducible $\Uehchi$-module  and $V$ is a quotient of $V'={\mathcal Ind}_T^{G,\chi}W$.  By \cite[Thm 5.1]{DCKP-regular}, we have $\dim V=\ell^{|\Phi^+|}=\dim V'$,  so $V={\mathcal Ind}_T^{G,\chi}W$. The general statement follows from Lemma~\ref{transitivity}.
\hfill$\Box$


\begin{proposition}
If $G$ is of type $A_{n-1}$ and $\mO_g^G$ is a subregular unipotent conjugacy class, then Conjecture~\ref{conj:gio} holds.
\end{proposition}
\pf The subregular unipotent conjugacy class $\mO_g^G$ in type $A_{n-1}$ is induced from the trivial class in a Levi subgroup of type $A_1$. It was shown in \cite[Thm. 3.11]{subregular} that all representations $V$ in $\spec(\Ueg)$ for which $\varphi(V)=\mO_g^G$ are induced from a representation of $\Uel$ where $\Uel$ is a reductive quantized enveloping algebra corresponding to a Levi subalgebra with semisimple part of type $A_1$.\hfill$\Box$

\begin{corollary}\label{cor:leq}
Conjecture~\ref{conj:gio} holds for $\mO_{su}^G$ whenever the derived subgroup of the Levi envelope $L$ of $G^{s}$ is of type $A_{a_1}\times\cdots\times A_{a_k}$ with $a_i\leq 2$. 
\end{corollary}
\pf By Theorem~\ref{thm:DCK-reduction} every irreducible module $V$ such that $\varphi(V)\in\mO_{su}^G$ is induced from an irreducible $\Uelchil$-module $L$. On the other hand, as in the proof of Theorem~\ref{thm:jordan} we see that every irreducible $\Uelchil$-module is, up to tensoring with a one-dimensional representation, a module $W$ for which $\varphi_L(W)$ is a unipotent class in $[L,\,L]$. By the discussion above, Conjecture~\ref{conj:gio} holds for $G=SL_2({\mathbb C})$ and $SL_3({\mathbb C})$ so we conclude by using Lemma~\ref{transitivity}.
\hfill$\Box$

\begin{corollary}Let $G=SL_n({\mathbb C})$ and let $s\in T$ be such that each eigenvalue of $s$ is repeated at most $3$ times. Then,
Conjecture~\ref{conj:gio} holds for $\mO_{su}^G$, for every $u\in G^{s}$. 
\end{corollary}
\pf Immediate from Corollary~\ref{cor:leq}\hfill$\Box$

\begin{remark}Most of the statements for $\Ueg$ are still valid or may be modified in order to hold  for $G$ not-necessarily simply-connected. In this case, the centralizer of a semisimple element should be replaced by its identity component. 
\end{remark}

\subsection{Acknowledgements}
 I am grateful to the organizers of the conference ``Hopf algebras and tensor categories'', where the content of this paper was illustrated,
 for  invitation, partial support and very warm hospitality.
Participation to the conference was also partially supported by Progetto di Ateneo CPDA105885 of the University of Padova. I wish to thank Prof. Istv\`an Heckenberger for his time, patience and the many useful discussions  on and around this topic during a visit to the University of Marburg in march 2012, supported by Italian PRIN 20097NBFW5-002, ``Azioni di gruppi: aspetti algebrici e geometrici''. Finally, I thank the referee for useful comments and remarks.


\begin{thebibliography}{00}


\bibitem{bo}{\sc W. Borho,}
\newblock{\em \"Uber Schichten halbeinfacher Lie-Algebren,}
\newblock Invent. Math., 65, 283--317 (1981/82). 

\bibitem{BK}{\sc W. Borho, H. Kraft,}
\newblock {\em \"Uber Bahnen und deren Deformationen bei linearen
  Aktionen reduktiver Gruppen,}
\newblock Comment. Math. Helvetici, 54, 61--104  (1979).








\bibitem{modp}{\sc N. Cantarini,}
\newblock{\em Mod-p reduction for quantum groups,}
\newblock J. Algebra, 202(1), 357--366 (1998).

\bibitem{spherical}{\sc N. Cantarini,}
\newblock{\em Spherical orbits and quantized enveloping algebras,}
\newblock Comm. Algebra, 27(7), 3439--3458 (1999). 

\bibitem{subregular}{\sc N. Cantarini,} 
\newblock{\em The quantized enveloping algebra $U_q(sl(n))$ at the roots of unity,}
\newblock Commun. Math. Phys., 211(1),  207--230 (2000). 


\bibitem{ccc}{\sc N.\ Cantarini, G.\ Carnovale, M.\ Costantini,}
\newblock {\em Spherical orbits and representations of $\Ueg$,}
\newblock Transformation\  Groups, 10(1), 29--62  (2005).


\bibitem{gio-espo}{\sc G.\ Carnovale, F. Esposito,}
\newblock{\em On sheets of conjugacy classes in good characteristic,}
\newblock IMRN, 2012(4), 810--828 (2012).






\bibitem{DCK1}{\sc C.\ De Concini, V.\ G.\ Kac,}
\newblock {\em Representations of quantum groups at roots of $1$,}
\newblock In:``Operator Algebras, Unitary Representations, Enveloping Algebras, and Invariant Theory'', Actes du colloque en l'honneur de Jacques Dixmier, Progress in Mathematics 92, Birkh\"auser, Boston, 471--506 (1990).


\bibitem{DCK}{\sc C.\ De Concini, V.\ G.\ Kac,}
\newblock {\em Representations of quantum groups at roots of $1$:
  reduction to the exceptional case,}
\newblock Infinite analysis, Part A, B (Kyoto,1991), 141--149,
Adv. Ser. Math. Phys., 16, World Sci. Publ., River Edge, NJ, (1992).  

\bibitem{DCKP}{\sc C.\ De Concini, V.\ G.\ Kac, C.\ Procesi,}
\newblock {\em  Quantum coadjoint action,}
\newblock J. Amer.\ Math.\ Soc.\ 5, 151--190 (1992).

\bibitem{DCKP-regular}{\sc C.\ De Concini, V.\ G.\ Kac, C.\ Procesi,}
\newblock {\em Some remarkable degenerations of quantum groups,}
\newblock Comm. Math. Phys. 157, 405--427 (1993).


\bibitem{DCM}{\sc C. De Concini, A. Maffei,}
\newblock{\em A generalized Steinberg section and branching rules for quantum groups at roots of $1$}
\newblock arxiv:1107.0248v1 (2011).

\bibitem{DCPRR}{\sc C. De Concini, C. Procesi, N. Reshetikhin, M. Rosso,}
\newblock{\em Hopf algebras with trace and representations,}
\newblock Invent. Math. 161(1), 1--44 (2005).


\bibitem{ela}{\sc  A. G. Elashvili,}
\newblock{\em Sheets of the exceptional Lie algebras,}
\newblock Issledovaniya po algebre, Tbilisi 171--194 Russian, (1985).










\bibitem{ke}{\sc G. Kempken,}
\newblock{\em Induced conjugacy classes in classical Lie-algebras,}
\newblock Abh. Math. Sem. Univ. Hamburg 53, 53--83 (1983).

\bibitem{kostant}{\sc B.\ Kostant}
\newblock{\em The principal three-dimensional subgroup and the Betti
  numbers of a complex simple Lie group,}
\newblock Amer. J. Math. 81, 973-1032 (1959). 

%






\bibitem{lusztig}{\sc G. Lusztig,}
\newblock{\em Intersection cohomology complexes on a reductive group,}
\newblock Invent. Math. 75(2), 205--272 (1984).

\bibitem{lusp}{\sc G. Lusztig, N. Spaltenstein,}
\newblock{\em Induced unipotent classes,}
\newblock J. London Math. Soc. 19(2), 41--52 (1979).

\bibitem{lu-braid}{\sc G. Lusztig,}
\newblock{\em On quantum groups,}
\newblock J. Algebra 131(2), 466--475, (1990).















\bibitem{spaltenstein}{\sc N. Spaltenstein,}
\newblock{\em Classes Unipotentes et Sous-groupes de Borel,}
\newblock LNM 946, Springer-Verlag (1982).



\bibitem{131}{\sc T.A.\ Springer, R.\ Steinberg,}
\newblock{\em Conjugacy classes,}
\newblock In: ``Seminar on algebraic groups and related finite groups''.
LNM 131, 167--266, Springer-Verlag, Berlin Heidelberg New York (1970).








\end{thebibliography}
\end{document}